\documentclass[12pt]{article}
\usepackage{amsthm}
\usepackage{amsmath}
\usepackage{amssymb}
\usepackage{epsfig}
\usepackage[latin5]{inputenc}

\numberwithin{equation}{section}
\newcommand{\gen}[1]{\partial_{#1}}

\DeclareMathOperator{\Sl}{sl} 
 
\DeclareMathOperator{\SL}{SL}

\newtheorem{prop}{Proposition}
\newtheorem{thm}{Theorem}

\newtheorem{defn}{Definition}

\begin{document}
\title{\bf
\Large A variable coefficient nonlinear Schr\"{o}dinger equation with a four-dimensional symmetry group and blow-up of its solutions}

\author{F. G\"{u}ng\"{o}r\thanks{Department of Mathematics, Faculty of Arts and Sciences, Do\u{g}u\c{s} University, 34722 Istanbul, Turkey, e-mail: fgungor@dogus.edu.tr and e-mail: mhasansoy@dogus.edu.tr},\and  M. Hasanov\footnotemark[1],\and
C. \"{O}zemir\thanks{Department of Mathematics, Faculty of Science and Letters,
Istanbul Technical University, 34469 Istanbul,
Turkey, e-mail: ozemir@itu.edu.tr}}


\date{\today}

\maketitle

\begin{abstract}
 A canonical variable coefficient nonlinear Schr\"{o}dinger equation with a four dimensional symmetry group containing $\SL(2,\mathbb{R})$ group as a subgroup is considered. This typical invariance is then used to transform by a symmetry transformation a  known solution  that can be derived by truncating its Painlev\'e expansion  and study  blow-ups of these solutions in the $L_p$-norm for $p>2$, $L_\infty$-norm and in the sense of
distributions.
\end{abstract}

\emph{Keywords:}   $\SL(2,\mathbb{R})$ invariance, variable coefficient nonlinear Schr\"{o}dinger equation, exact solutions, blow-up

\emph{AMS subject Classifications:} Primary 35Q55, 35B44, 35B06; Secondary 35A25

\section{Introduction}

It is well known that the linear heat and Schr\"{o}dinger equations for $u\in \mathbb{R}$, $\psi\in \mathbb{C}$
 $$ u_t=u_{xx},\quad i\psi_t+\psi_{xx}=0,\quad x\in \mathbb{R},\quad t>0$$
 have isomorphic Lie symmetry groups.  The symmetry group with the infinite-dimensional ideal reflecting linearity factored out can be written as a semidirect product of the three-dimensional Heisenberg group $\mathsf{H}$ and $\SL(2,\mathbb{R})$
$$ G=\mathsf{H}\ltimes \SL(2,\mathbb{R}).$$ On the other hand, among the modular class of nonlinear Schr\"{o}dinger equations in one space dimension
\begin{eqnarray}\label{modular}
i\psi_t+\psi_{xx}=F(|\psi|)\psi,
\end{eqnarray}
the only one preserving the symmetry group of the linear Schr\"{o}dinger equation  except with infinite-dimensional symmetry of its linear counterpart is the quintic Schr\"{o}dinger equation
\begin{equation}
i\psi_t+\psi_{xx}=\lambda|\psi|^4\psi.
\end{equation}

In this article we look at a variable coefficient extension of the one-dimensional cubic NLSE (nonlinear Schr\"{o}dinger equation)
\begin{eqnarray}\label{VCNLS}
\begin{split}
& i\psi_t+\psi_{xx}+g(x,t) |\psi|^2 \psi+h(x,t) \psi=0,\\
& g=g_1+{i}g_2, \quad h=h_1+{i}h_2,\quad
g_j,h_j\in\mathbb{R}, \quad j=1,2, \quad g_1\neq 0.
\end{split}
\end{eqnarray}
Variable coefficient extensions of nonlinear evolution type
equations tend to arise in  cases when less idealized conditions
such as inhomogeneities and variable topographies are assumed in
their derivation. The reader is referred to \cite{Ozemir10} and the references therein for several physically motivated applications.
Symmetry classification of \eqref{CNLS} was given in \cite{Gagnon93}. We mention that the dimension of the maximal symmetry group is $\dim G=5$ and is achieved only when the coefficients are constant and additionally $h=0$ which is reduced to nothing more than the usual NLS equation. The symmetry group is isomorphic to the group of one-dimensional  extended Galilei similitude algebra. We are particularly interested in the case when at least the $\SL(2,\mathbb{R})$  invariance is contained in the full symmetry group of \eqref{VCNLS}. This is usually referred to  as the pseudo-conformal invariance in the context of qualitative analysis of PDEs.

We emphasize that this  invariance manifests itself as a subgroup of the full symmetry groups of the NLSE, Davey-Stewartson (DS) equations and their
possible generalizations \cite{Gungor06} and has been successfully applied to  investigate blowup formation in
these nonlinear evolution models \cite{Cipolatti01, Ozawa92, Cazenave91, Eden06}.  A study of self-similar solutions of the pseudo-conformally invariant nonlinear Schr\"{o}dinger equation can be found in \cite{Kavian94}.  The point is that when the variable coefficients  are allowed
in these equations, this typical symmetry is mostly destroyed.   Our intention here is  to detect the subcases in which such
a symmetry remains intact in variable coefficients variants of the NLS equations. This will make it possible to generate new nontrivial solutions from known ones.

We quote the following result from \cite{Gagnon93} and note that
throughout the paper, any solution will be understood as a
pointwise solution in the classical sense.
\begin{prop}
Any equation of the form \eqref{VCNLS} containing  $\SL(2,\mathbb{R})$ symmetry group as a subgroup can be transformed by point transformations to the canonical form
\begin{equation}\label{CNLS}
i\psi_t+\psi_{xx}+(\epsilon+i \gamma)\frac{1}{x} |\psi|^2 \psi+(h_1+i h_2)\frac{1}{x^2} \psi=0,\quad x\in \mathbb{R}\setminus\{0\},
\end{equation}
where $\epsilon=\pm 1$, $\gamma$, $h_1$ and $h_2$ are arbitrary real constants.
\end{prop}

Note that  like all equations in the modular class \eqref{modular}, equation \eqref{VCNLS} and also \eqref{CNLS} are always invariant under the constant change of phase of $\psi$ (gauge-invariance) while leaving the $(x,t)$ coordinates unchanged.
We represent the phase and modulus of $\psi$ by $\rho$, $\omega$ writing $\psi=\rho(x,t)\exp{(i\omega(x,t))}$.

\begin{prop}
The symmetry algebra $L$ of \eqref{CNLS} is four-dimensional and spanned by the vector fields
$$ T=\gen t,  \quad D=2t\gen t+x\gen x-\frac{1}{2} \rho \gen \rho,\quad C=t^2\gen t+xt\gen x-\frac{1}{2}t\rho\gen \rho+\frac{x^2}{4}\gen \omega,\quad W=\gen \omega.$$
The commutators among $T, D, C$ satisfy
$$ [T,D]=2T, \quad [T,C]=D, \quad [D,C]=2C$$ and $W$ is the center element, namely commutes with all other elements.
The Lie algebra  $L$ has the direct sum structure
$$ L=\Sl(2,\mathbb{R})\oplus \mathbb{R}.$$
\end{prop}
The elements $T$, $D$, $C$, $W$ generate time translations, scaling, (pseudo)-conformal
and gauge transformations, respectively. A significant consequence of these transformations  is the group action on the solutions given by the following.

\begin{prop}
If $\psi_0(x,t)$ is a solution of \eqref{CNLS}, then so is
$$\psi(x,t)=(a+bt)^{-1/2}\, \mathrm{e}^{{i}\frac{bx^2}{4(a+bt)}}
\,\psi_0\left(\frac{x}{a+bt},\frac{c+dt}{a+bt}\right)$$
for $ad-bc=1$.
\end{prop}
\proof{By exponentiating the infinitesimal
generators $T,D,C$ (i.e. solving  Cauchy problems for these vector fields) and then composing the corresponding group transformations  we find the above
$\SL(2,\mathbb{R})$ action on the solutions. Note that the action  corresponding to $t$ generates M\"{o}bius transformations of $t$.}

Based on this result, our main purpose is to transform one known
solution   of the original equation to a more complicated one by the transformations of  the $\SL(2,\mathbb{R})$ symmetry
and then  choose group parameters $a,b,c$ appropriately and pass to limit   of
the wave function $\psi$ as $t\to T^{-}$ for some finite time $T$ in the $L_p$-norm and distributional sense as well.

We now use a truncation  approach to obtain a special exact solution of \eqref{CNLS}.
We are going to find this special explicit solution   by truncating
its Painlevé expansion at the first term.  For convenience we write \eqref{CNLS}
together with its complex conjugate as the system
\begin{eqnarray}\label{sys}
\begin{split}
iu_t+u_{xx}+(\epsilon+i \gamma)\frac{1}{x} u^2 v+(h_1+i h_2)\frac{1}{x^2} u=0,\\
-iv_t+v_{xx}+(\epsilon-i \gamma)\frac{1}{x} u v^2+(h_1-i h_2) \frac{1}{x^2} v=0.\\
\end{split}
\end{eqnarray}
Here $u$ was employed instead of $\psi$ and $v$ denotes its
complex conjugate, but in this setting they are viewed as
independent functions. We first show that \eqref{sys} does not
pass the Painlevé test for PDEs and then proceed to obtain an
exact solution afterwards.

A partial differential equation is said to have the Painlevé
property if all its solutions are single valued around any
non-characteristic movable singularity manifold. If this
singularity manifold is denoted by $\Phi(x,t)=0$ (actually a curve
in this case), we shall look for solutions of the system
~\eqref{sys} in the form of a Laurent expansion  and  we expand
\begin{equation}\label{exp}
u(x,t)=\sum_{j=0}^{\infty} \; u_j(x,t)\Phi^{\alpha+j}(x,t), \quad
v(x,t)=\sum_{j=0}^{\infty} \; v_j(x,t)\Phi^{\beta+j}(x,t),
\end{equation}
where  $u_0, v_0\ne 0$ and $u_j, v_j, \Phi(x,t)$ are analytic
functions. $\alpha$ and $\beta$ are negative integers to be
determined from the leading order analysis, so as to ensure
absence of essential singularities and branch points in all
solutions. For the determination of leading orders $\alpha$ and
$\beta$, we substitute $u\sim u_0 \Phi^\alpha$ and $v\sim
v_0\Phi^\beta$ in \eqref{sys} and see that by balancing the terms
of smallest order
\begin{equation}\label{ab}
\alpha+\beta=-2\\
\end{equation}
and
\begin{equation}\label{u0v0}
u_0v_0=-\alpha(\alpha-1)\frac{x}{\epsilon+i\gamma}\Phi_x^2=
-\beta(\beta-1)\frac{x}{\epsilon-i\gamma}\Phi_x^2
\end{equation} must hold. \eqref{ab} allows the negative integers
$\alpha=-1$ and $\beta=-1$ and with these leading orders,
\eqref{u0v0} forces $\gamma=0$. After determination of the
leading orders, we substitute \eqref{exp} into \eqref{sys}.
Equating to zero the coefficient of $\Phi^{-3+j}, \,j\geq1,
j\in\mathbb{N}$, we arrive at a linear system
\begin{equation}\label{det}
Q(j)\left( \begin{array}{c}
        u_j\\
        v_j
       \end{array}\right)=
 \left( \begin{array}{c}
        F_j\\
        G_j
       \end{array}\right)
\end{equation}
from which the coefficients $u_j, \, v_j$ can be found. Those
values of indices $j$ for which $\det Q(j)=0$ are called \emph{resonances}. In
order that the expansion \eqref{exp} includes correct number of
arbitrary functions as required by the Cauchy-Kovalewski theorem
(where $\Phi(x,t)$ should be one of the arbitrary functions), some
consistency conditions  at resonances
must be satisfied. In the general case \eqref{VCNLS}, these constraints force the coefficients to be properly related (See \cite{Ozemir10} for details).  For \eqref{CNLS}, this is not the case, namely it cannot pass the
Painlevé test. It thus fails to satisfy the necessary condition for
the equation to have Painlevé property.

Notwithstanding this fact, application of the Painlevé expansion to nonintegrable PDEs like \eqref{CNLS} (or more properly partial integrable) can allow particular explicit solutions
 to be obtained by truncating the expansion. This approach  imposes constraints on the arbitrary functions
 and the function $\Phi$ as a result of  compatibility of an overdetermined
PDE system.

Before truncating the series \eqref{exp} at some index $j=N$, we
first weaken the condition that $\alpha$ and $\beta$ assumes only
integer values. When we solve \eqref{ab} and \eqref{u0v0}
together, we find that for $\gamma \neq 0$
\begin{equation}\label{ab1}
\alpha=-1-i\delta, \quad \beta=-1+i\delta;
\qquad \delta=\frac{-3\epsilon\pm\sqrt{8\gamma^2+9}}{2\gamma}
\end{equation}
and \eqref{u0v0} is equivalent to
\begin{equation}\label{u0v01}
u_0v_0=-\frac{3\delta}{\gamma}\,x\,\Phi_x^2.
\end{equation}
Truncating the Painlevé expansion at the first term ($j=0$), we
assume a solution of the form
\begin{equation}
u(x,t)=u_0(x,t)\Phi(x,t)^{-1-i\delta}, \quad
v(x,t)=v_0(x,t)\Phi(x,t)^{-1+i\delta}.
\end{equation}
We substitute these in \eqref{sys} and set the coefficients of the
terms $\Phi^{-3\pm i\delta}, \Phi^{-2\pm
i\delta}, \Phi^{-1\pm i\delta}$  equal to zero. The
condition at the order $\Phi^{-3\pm i\delta}$ is
equivalent to \eqref{u0v01} and terms of order $\Phi^{-2\pm
i\delta}, \Phi^{-1\pm i\delta}$ disappear if $\Phi$, $u_0$, $v_0$ are chosen to be
\begin{align}
\Phi(x,t)&=\left(\frac{x}{k_4t+k_1}\right)^{2/3}+k_2\\
u_0(x,t)&=A_1\frac{x^{1/6}}{(k_4t+k_1)^{2/3}}\exp\left[i\Bigl(\frac{k_4x^2}{4(k_4t+k_1)}+k_3\Bigr)\right]\\
v_0(x,t)&=A_2\frac{x^{1/6}}{(k_4t+k_1)^{2/3}}\exp\left[-i\Bigl(\frac{k_4x^2}{4(k_4t+k_1)}+k_3\Bigr)\right]
\end{align}
for arbitrary real constants $A_1,A_2$ and $k_1,...,k_4$ with the constraints
$h_1=5/36$ and $h_2=0$. We require to have $v=u^*$, which implies
$A_1=A_2$. Let us rename this constant $A$. When we check the
condition \eqref{u0v01}, we find that
$A^2=\frac{-4\delta}{3\gamma}$. Since this square must be
positive, it is necessary that
\begin{equation}
\frac{\delta}{\gamma}=\frac{-3\epsilon\pm\sqrt{8\gamma^2+9}}{2\gamma^2}<0.
\end{equation}
For both $\epsilon=\pm1$, it is seen that the minus sign must be
picked in the formula  for $\delta$ of \eqref{ab1}. Having
found a consistent truncation, we can write the solution to
\eqref{CNLS} as
\begin{equation}\label{td}
u(x,t)=\frac{Ax^{1/6}}{x^{2/3}+k_2(k_4t+k_1)^{2/3}}
\exp\left[i\left(\frac{k_4x^2}{4(k_4t+k_1)}
-\delta\ln\left(\frac{x^{2/3}}{(k_4t+k_1)^{2/3}}+k_2\right)+k_3\right)\right].
\end{equation}
Choosing $k_4=0$, we obtain the stationary solution
\begin{eqnarray}\label{stat}
u(x)=\frac{Ax^{1/6}}{x^{2/3}+k_1^{2/3} k_2}
\exp\left[i\left(
-\delta\ln(x^{2/3}+k_1^{2/3}k_2)+k_3\right)\right],
\end{eqnarray}
where
\begin{equation}\label{constraint}
h_1=\frac{5}{36}, \quad h_2=0, \quad
A=(-\frac{4\delta}{3\gamma})^{1/2}, \quad
\delta=\frac{-3\epsilon-\sqrt{8\gamma^2+9}}{2\gamma},
\end{equation}
and $k_1, k_2, k_3$ ($k_3$ is relabelled) are arbitrary real constants.

\noindent We summarize:
\begin{prop}
The following solves equation \eqref{CNLS} for arbitrary constants $k_1,
k_2, k_3$ and for the parameters $h_1=5/36$, $h_2=0$
\begin{equation}
u(x,t)=A\,\frac{x^{1/6}}{x^{2/3}+k_1^{2/3} k_2}
\exp\left[i\left(
-\delta\ln(x^{2/3}+k_1^{2/3}k_2)+k_3\right)\right],
\end{equation}
where $A=(-\frac{4\delta}{3\gamma})^{1/2}$ and
$\delta=\frac{-3\epsilon-\sqrt{8\gamma^2+9}}{2\gamma}$.
\end{prop}

\section{Transforming solutions  by $\SL(2,\mathbb{R})$ group}
\subsection*{Blowup in the $L_p$, $L_\infty$ norms and in the distributional sense}

Now we would like to illustrate how the
$\SL(2,\mathbb{R})$ group action can be useful in establishing blow-up
profiles of initial value problems for variable coefficient NLS
equations just as they were used for their constant coefficient
counterparts.

Let $u(x)$ be the stationary solution to \eqref{CNLS}, defined by
\eqref{stat}-\eqref{constraint}. We set $\psi_0(x):=u(x)$ and use Proposition 3. By
this proposition for arbitrary $a, b\in \mathbb{R}$
\begin{equation}\label{2.1}
\psi(x,t)=(a+bt)^{-1/2}\,
\exp{\Bigl({\frac{ibx^2}{4(a+bt)}}\Bigr)}
\,\psi_0\left(\frac{x}{a+bt}\right)
\end{equation}
is also a solution to equation \eqref{CNLS}.

We can assume $a>0, b<0$ and
denote $\varepsilon:=a+bt=b(t- \frac{-a}{b})=b(t-T),$ where
$T=-\frac{a}{b}>0$. Hence, $t\rightarrow T^-\Leftrightarrow
\varepsilon \rightarrow 0^+.$ By using this notation we can write
solution \eqref{2.1} in the form
$$\psi_{\varepsilon}(x)=\varepsilon^{-1/2}\exp{\Bigl(\frac{ibx^2}{4\varepsilon}\Bigr)}
\psi_0\bigl(\frac{x}{\varepsilon}\bigr).$$

We are going to show that these solutions will blowup in the
$L_p$-norm when $p>2$, $L_\infty$-norm and in the sense of
generalized functions, respectively.

\paragraph{\bf Note:} In the following theorems, we do not impose an initial condition but rather we limit to the one  as dictated by the solution \eqref{2.1} in the form
$$\psi(x,0)=\frac{1}{\sqrt{a}}\exp{\Bigl(\frac{ibx^2}{4a}\Bigr)}\psi_0\bigl(\frac{x}{a}\bigr).$$ Also, for a blow-up at some finite time we have to fix $b$ figuring in \eqref{2.1}.

\paragraph{$L_p$-blow-up solutions.}

\begin{thm}
 For any $T>0$ there is a solution $\psi(x,t)$ to equation \eqref{CNLS} such
that
\begin{equation}
\lim_{t\rightarrow
T^-}\|\psi(x,t)\|_p=+\infty\quad\hbox{for\,\,all}\,\,p>2,
\end{equation}
where $\|\psi(x,t)\|_p=
\Bigl(\int_{-\infty}^{+\infty}|\psi|^p\,dx\Bigr)^{1/p}$.
\end{thm}
\proof  Let $T>0$ be a  finite
time. We can always arrange two numbers $a>0$ and $b<0$ such that
$T=-\frac{a}{b}$. Setting  these numbers in \eqref{2.1} we get the
function $\psi(x,t)$ which will be instrumental in the proof. We
rewrite this function, as given above, in the form

$$\psi_{\varepsilon}(x)=\varepsilon^{-1/2}\exp{\Bigl(\frac{ibx^2}{4\varepsilon}\Bigr)}\psi_0\bigl(\frac{x}{\varepsilon}\bigr).$$
Then we have
\begin{equation}
\lim_{\varepsilon\rightarrow
0^+}\|\psi_\varepsilon(x)\|_p=\lim_{t\rightarrow T^-}\|\psi
(x,t)\|_p.
\end{equation}
 By the definition of
$\psi_0(x)$ we have
$$
\psi_0\bigl(\frac{x}{\varepsilon}\bigr)
=\frac{A\varepsilon^{-\frac{1}{6}}x^{\frac{1}{6}}}{\varepsilon^{-\frac{2}{3}}x^{\frac{2}{3}}+
C} \exp\bigl(i\varphi (x,\varepsilon)\bigr),
$$
where $ \varphi(x,\varepsilon)= -\delta
\ln\Bigl({\bigl(\frac{x}{\varepsilon}}\bigr)^{\frac{2}{3}}+k_1^{\frac{2}{3}}k_2+k_3\Bigr)$,
 $C= k_1^{\frac{2}{3}}k_2$ and we choose $C>0$. Thus,
 $$
\psi_\varepsilon(x) =\frac{Ax^{\frac{1}{6}}}{x^{\frac{2}{3}}+
\varepsilon^{\frac{2}{3}}C}
\exp\bigl[i\Bigl(\frac{bx^2}{4\varepsilon}+\varphi
(x,\varepsilon)\Bigr)\bigr],
 $$
and
$$
\int_{-\infty}^{\infty}|\psi_\varepsilon(x)|^pdx=2\int_0^{\infty}\frac{A^px^{p/6}}{{\bigl(x^{2/3}+\varepsilon^{2/3}C}\bigr)^p}\,dx.
$$
By (2.3) $\lim_{t\rightarrow T^-}\|\psi(x,t)\|_p=+\infty$ if and only if:
\begin{enumerate}

\item[ i)] $\int_0^{\infty}|\psi_\varepsilon(x)|^p$ is finite for all
$p>2$ and $\varepsilon >0$,

\item[ ii)] $\lim_{\varepsilon\rightarrow 0^+}
\int_0^{\infty}|\psi_\varepsilon(x)|^pdx= \infty.$
\end{enumerate}
Let us substitute $x=\varepsilon y$. Then
\begin{equation}
\begin{split}
\int_{-\infty}^{\infty}|\psi_\varepsilon(x)|^pdx
&=
2A^p\int_0^{\infty}\frac{\varepsilon^{\frac{p}{6}}y^{\frac{p}{6}}}{{\varepsilon^{\frac{2p}{3}}\bigl(y^{\frac{2}{3}}+C}\bigr)^p}\varepsilon dy=
\frac{2A^p}{\varepsilon^{(p-2)/2}}\int_0^{\infty}\frac{y^{\frac{p}{6}}}{{\bigl(y^{\frac{2}{3}}+C}\bigr)^p}dy=\\
&=\frac{2A^p}{\varepsilon^{(p-2)/2}}\Bigl[\int_0^1\frac{y^{\frac{p}{6}}}{{\bigl(y^{\frac{2}{3}}+C}\bigr)^p}dy+
\int_1^{+\infty}\frac{y^{\frac{p}{6}}}{{\bigl(y^{\frac{2}{3}}+C}\bigr)^p}dy\Bigr],
\end{split}
\end{equation}
where
$\int_0^1\frac{y^{\frac{p}{6}}}{{\bigl(y^{\frac{2}{3}}+C}\bigr)^p}dy$
is convergent by continuity of the function. On the other hand,
$$
\frac{y^{\frac{p}{6}}}{{\bigl(y^{\frac{2}{3}}+C}\bigr)^p}\leq
\frac{y^{\frac{p}{6}}}{ y^{\frac{2p}{3}}}=\frac{1}{
y^{\frac{p}{2}}}.
$$
Consequently,
$\int_1^{\infty}\frac{y^{\frac{p}{6}}}{{\bigl(y^{\frac{2}{3}}+C}\bigr)^p}dy\leq\int_1^{\infty}\frac{1}{
y^{\frac{p}{2}}}dy$ implies that the integral
$\int_1^{\infty}\frac{y^{\frac{p}{6}}}{{\bigl(y^{\frac{2}{3}}+C}\bigr)^p}dy$
is convergent for all $p>2$. Hence,
$\int_0^{\infty}|\psi_\varepsilon(x)|^pdx$ is convergent for
$p>2$. Finally, by taking the limit we find
$$
\lim_{\varepsilon\rightarrow
0^+}\int_{-\infty}^{\infty}|\psi_\varepsilon(x)|^pdx=\lim_{\varepsilon\rightarrow
0^+}
\frac{2A^p}{\varepsilon^{(p-2)/2}}\int_0^{\infty}\frac{y^{\frac{p}{6}}}{{\bigl(y^{\frac{2}{3}}+C}\bigr)^p}dy=
\infty.
$$
\paragraph{Blow-up  solutions in  $L_\infty$-\,norms.}

In the following theorem we prove that the above defined
solutions $\psi_\varepsilon(x)$ will blowup in $L_\infty$- norm
too.

\begin{thm}
For any $T>0$ there is a solution $\psi(x,t)$ to equation \eqref{CNLS}
such that
\begin{equation}
\lim_{t\rightarrow T^-}\|\psi(x,t)\|_\infty=\infty,
\end{equation}
where $\|\psi(x,t)\|_\infty= {\rm ess}\sup_{x\in
[0,\infty)}|\psi(x,t|,\,\,t<T$.
\end{thm}

\proof Let again $\psi (x,t)$ be the function defined by \eqref{2.1} and
$\psi_\varepsilon(x)$ be the above defined function. By the
construction of $\psi_\varepsilon(x)$
$$
\lim_{\varepsilon\rightarrow
0^+}\|\psi_\varepsilon(x)\|_\infty=\lim_{t\rightarrow
T^-}\|\psi(x,t)\|_\infty
$$
and therefore we are done if we can show that $\lim_{\varepsilon\rightarrow
0^+}\|\psi_\varepsilon(x)\|_\infty=\infty$. We have
$|\psi_\varepsilon(x)|=\frac{A|x|^{\frac{1}{6}}}{x^{\frac{2}{3}}+
\varepsilon^{\frac{2}{3}}C}$. Let $A=C=1$. $|\psi_\varepsilon(x)|$ is an even function so that we can restrict ourselves to the interval $[0,\infty)$. Since
$|\psi_\varepsilon(x)|$ is continuous on  $[0,\infty)$,
$|\psi_\varepsilon(0)|=0$ and
$\lim_{x\rightarrow\infty}|\psi_\varepsilon(x)|=0$, then there
exists $x_0\in (0,\infty)$ such that
$$
\|\psi_\varepsilon(x)\|_\infty= {\rm ess}\sup_{x\in
[0,\infty)}|\psi_\varepsilon(x|=\max_{x\in
[0,\infty)}|\psi_\varepsilon(x)|=\frac{x_0^{\frac{1}{6}}}{x_0^{\frac{2}{3}}+
\varepsilon^{\frac{2}{3}}C}.
$$
A simple computation yields $x_0=\frac{\varepsilon}{\sqrt{27}}$.
Therefore,
$$
\|\psi_\varepsilon(x)\|_\infty=
C\frac{\varepsilon^{\frac{1}{6}}}{{\varepsilon^\frac{2}{3}}}=
\frac{C}{\sqrt{\varepsilon}}.
$$
Consequently,
$$
\lim_{\varepsilon\rightarrow
0^+}\|\psi_\varepsilon(x)\|_\infty=\infty.
$$

\paragraph{$\delta$-Blowup solutions in the sense of generalized
functions.}

Let $D= C_0^\infty(0,\infty)$ be  the space of infinitely
differentiable functions with  compact support in $(0,\infty)$.
The dual of $D$ is called the space of generalized functions and
is denoted by  $D'$.

\begin{defn}\label{def1}
Let $f,\,f_n\in D'$. We say that the sequence $f_n$ converges to
$f$ if and only if $\langle f_n,\varphi\rangle \rightarrow\langle
f,\varphi\rangle$ for all $\varphi\in D$, where $\langle
f,\varphi\rangle$ denotes the value of the functional $f$ at
$\varphi$.
\end{defn}
Now we present  $\delta$-blowup solutions in the sense of
generalized functions.
\begin{thm}
For all $p>2$
$$
\varepsilon^{(p-2)/2}|\psi_\varepsilon(x)|^p\rightarrow
K\delta(x) \quad  \hbox{as}\,\, \varepsilon\rightarrow 0^+,
$$
where $K=
A^p\int_{-\infty}^{\infty}\frac{|y|^{\frac{p}{6}}}{{\bigl(y^{\frac{2}{3}}+C}\bigr)^p}dy,\,\,C>0$
and $\delta(x)$ denotes the Dirac distribution at the origin.
\end{thm}

\proof By Definition \ref{def1} we have to show that
$$
\lim_{\varepsilon\rightarrow
0^+}\int_{-\infty}^{\infty}\varepsilon^{(p-2)/2}|\psi_\varepsilon(x)|^p\varphi(x)dx=
K\varphi(0)\quad \text{for all} \; \varphi\in D.
$$
Evidently,
$$
\varepsilon^{(p-2)/2}\int_{-\infty}^\infty|\psi_\varepsilon(x)|^p\varphi(x)\,dx=
\varepsilon^{(p-2)/2}\int_{-\infty}^\infty\frac{A^p|x|^{p/6}}{{\bigl(x^{2/3}+\varepsilon^{2/3}C}\bigr)^p}\varphi(x)\,dx.
$$
By setting the substitution $x=\varepsilon y$ we obtain
\begin{equation}
\begin{split}
&\varepsilon^{(p-2)/2}\int_{-\infty}^\infty|\psi_\varepsilon(x)|^p\varphi(x)\,dx=
\varepsilon^{(p-2)/2}\int_{-\infty}^\infty\frac{\varepsilon^{\frac{p}{6}}A^p|y|^{p/6}}{\varepsilon^{\frac{2p}{3}}{\bigl(y^{2/3}+C}\bigr)^p}\varphi(\varepsilon y)\,\varepsilon dy\\
&=\varepsilon^{(p-2)/2}\int_{-\infty}^\infty\frac{A^p|y|^{p/6}}{{\varepsilon^{(p-2)/2}\bigl(y^{2/3}+C}\bigr)^p}\varphi(\varepsilon
y)\,dy=
\int_{-\infty}^\infty\frac{A^p|y|^{p/6}}{{\bigl(y^{2/3}+C}\bigr)^p}\varphi(\varepsilon
y)\,dy .
\end{split}
\end{equation}
Let $f_\varepsilon(y):=
\frac{A^p|y|^{p/6}}{{\bigl(y^{2/3}+C}\bigr)^p}\varphi(\varepsilon
y)$. The sequence $f_\varepsilon(y)$ satisfies the following two
conditions:

i) $|f_\varepsilon(y)|\leq C_\varepsilon
\frac{A^p|y|^{p/6}}{{\bigl(y^{2/3}+C}\bigr)^p}\in
L_1(-\infty,\infty)$ if $p>2$,

ii) $\lim_{\varepsilon\rightarrow 0^+}f_\varepsilon(y)=
\frac{A^p|y|^{p/6}}{{\bigl(y^{2/3}+C}\bigr)^p}\varphi(0)$.

\noindent Then by Lebesgue's dominated convergence theorem we obtain
$$
\lim_{\varepsilon\rightarrow 0^+}\int_{-\infty}^\infty
f_\varepsilon(y)dy=
\varphi(0)\int_{-\infty}^\infty\frac{A^p|y|^{p/6}}{{\bigl(y^{2/3}+C}\bigr)^p}\,dy=
K\varphi(0).
$$
By Definition 1, this implies that
$$
\varepsilon^{(p-2)/2}|\psi_\varepsilon(x)|^p\rightarrow K\delta(x),
$$
as $\varepsilon\rightarrow 0^+$.

\noindent{\bf Remark:} The above argument indicates that the singular behavior of the solution is even much worse than in the usual distributional sense as was  done in Ref. \cite{Ozawa92}.


\end{document}